\documentclass[10pt]{article}
 \usepackage{latexsym}
 \usepackage{amsmath}
 \usepackage{amsfonts}

         \newtheorem{theorem}{Theorem}[section]

         \newcommand{\solidbox}{\vrule width.6em height.5em
 depth.1em\relax}

 \newcommand{\qed}{\unskip\enskip\null\nobreak\hfill\solidbox
 \par}
         
         \newcommand{\N}{\mathbb{N}}
         \newcommand{\R}{\mathbb{R}}
         \renewcommand{\chi}{\bf{1}}

     \newcommand{\C}{\mathbb{C}}
     \renewcommand{\H}{\mathbb{H}}

\begin{document}

\title{Remarks on Naimark's duality}
\author{Wojciech Czaja\\ Institute of Mathematics, University of Wroc{\l}aw \\ Pl. Grunwaldzki 2/4, 50-384 Wroc{\l}aw, Poland\\and\\
Department of Mathematics, University of Vienna\\ Nordbergstrasse 15, 1090 Wien, Austria}
\footnotetext{\textit{Math Subject Classifications.} 42C40}
\footnotetext{The author is supported by Marie Curie Intra-European Fellowship FP6-2002-500685}

\maketitle
\vskip-1.5cm

\begin{abstract}
We present an extension of Naimark's duality principle which states that complete systems in a Hilbert space are projections of $\omega$-linearly independent systems of elements of an ambient Hilbert space. This result is presented in the context of other known extensions of Naimark's theorem.
\end{abstract}


\section{Introduction}

In 1940 M.~A.~Naimark \cite{N} proved the following result (stated here in a more modern language).

\begin{theorem}
\label{tn}
Any tight frame in a Hilbert space $\H$ is an orthogonal projection of an orthonormal basis of an ambient Hilbert space $\H'$, $\H \subset \H'$.
\end{theorem}

This result was later extended independently by Casazza, Han, and Larson, \cite{CHL} Theorem 7.6, and by Kashin and Kulikova, \cite{KK} Theorem 1, to hold for arbitrary frames and Riesz bases. In fact, Theorem 7.6 in \cite{CHL} shows that the property of a frame being a projection of a Riesz basis is directly connected to the notions of excesses and deficits of frames, see \cite{BCHL1} and \cite{BCHL2} for more on this subject.

In his recent work, Terekhin \cite{T} showed a further extension of the result of Naimark, which holds for representation systems and for Schauder bases. Terekhin's work is also more general in the sense that his results are stated for Banach spaces rather than Hilbert spaces.

In the main result of this paper, Theorem \ref{t1}, we present the above results combined with our own extensions. This way we obtain yet another duality principle in the theory of representation systems. This principle was motivated by the duality theory which holds for Gabor systems, i.e., systems generated by time-frequency shifts of a single function, see, e.g., \cite{RS}. We also refer the interested reader to a recent survey of analogous results in frame theory \cite{CKL}.

\section{Definitions and Main Result}

Let $\H$ be a separable, infinite-dimensional Hilbert space over $\C$. We say that a collection $\{ f_k : k=1, \ldots\}\subset \H$ of vectors is a {\it frame} for $\H$, with the {\it frame bounds} $A$ and $B$, if
\[
\forall\; f \in \H, \quad
A\| f \|_\H^2 \le \sum_{k\in\N} |\langle f,f_k \rangle|^2 \le B\| f \|_\H^2.
\]
We say that $A$ and $B$ are the lower and the upper frame bounds, respectively.
A frame is {\it tight} if $A=B$. A system that satisfies only the upper inequality in the above estimate is called a {\it Bessel system}. We say that a frame is {\it exact} if it is no longer a frame after removal of any of its elements. Exact frames are also called {\it Riesz bases}. An equivalent characterization of Riesz bases is commonly used, see, e.g., \cite{C}. It states that a collection $\{f_k : k=1, \ldots\}\subset \H$ of vectors is a Riesz basis if it is complete (i.e., $\overline{\text{\rm span}\; \{f_k : k=1, \ldots\}} = \H$) and there exist constants $0< A\le B<\infty$ such that
\begin{equation}
\label{upRiesz}
\forall\; c\in l^2(\N), \quad A \|c\|_2^2 \le \Big\| \sum_{k\in \N} c(k) f_k \Big\|_2^2 \le B\|c\|_2^2.
\end{equation}
When only the second inequality in the above holds, we say that the system $\{f_k : k=1, \ldots\}\subset \H$ satisfies the {\it upper Riesz inequality}. Such collections are in particular Bessel systems, e.g., Exercise 3.9 in \cite{C}.

A {\it representation system} $\{ f_k : k=1, \ldots\}\subset \H$ is a system such that for any $f \in \H$ there exists a sequence
$c = \{c(k): k = 1, \ldots\} \subset \C$ such that the series
\begin{equation}
\label{series}
\sum_{k \in \N} c(k) f_k
\end{equation}
converges to $f$ in the norm of $\H$. The space of sequences for which the above series converges is called the {\it coefficient space} of the representation system $\{f_k : k=1, \ldots\}$. A representation system with the property that the for each $f \in \H$ there exists a unique number sequence $c = \{c(k): k = 1, \ldots\}$ such that (\ref{series}) converges to $f$ in $\H$, is called a {\it Schauder basis}.

We define the {\it coefficient space of zero-series} to be the closed subspace $N$ of the coefficient space, consisting of these sequences for which the expansion (\ref{series}) converges to $0 \in \H$.
Finally, a collection $\{ f_k : k=1, \ldots\}\subset \H$ is $\omega$-linearly independent for $l^2(\N)$ if the fact that there exists $c \in l^2(\N)$ such that $\sum_{k\in\N} c_kf_k =0$ implies that $c = 0$.

\begin{theorem}
\label{t1}
Let $\H$ be a separable Hilbert space.

$\boldsymbol{a.}$ The collection $\{f_j:j \in \N\}\subset \H$ is a tight frame for $\H$ if and only if there exists a Hilbert space $\H'$ containing $\H$ and an orthogonal basis $\{e_j:j \in \N\}\subset \H'$ for $\H'$ such that
\[
\forall\; j\in \N, \quad P(e_j) = f_j.
\]

$\boldsymbol{b.}$ The collection $\{f_j:j \in \N\}\subset \H$ is a frame for $\H$ with frame constants $A \le B$ if and only if there exists a Hilbert space $\H'$ containing $\H$ and a Riesz basis $\{e_j:j \in \N\}\subset \H'$ for $\H'$ with constants the same constants $A\le B$ such that
\[
\forall\; j\in \N, \quad P(e_j) = f_j.
\]

$\boldsymbol{c.}$ Assume that the space of zero-series for $\{f_j:j \in \N\}\subset \H$ is complemented in the coefficient space of $\{f_j:j \in \N\}$. Then, $\{f_j:j \in \N\}$ is representation system in $\H$ if and only if there exists a Hilbert space $\H'$ containing $\H$ and a Schauder basis $\{e_j:j \in \N\}\subset \H'$ for $\H'$ such that
\[
\forall\; j\in \N, \quad P(e_j) = f_j.
\]

$\boldsymbol{d.}$ Assume that the collection $\{f_j:j \in \N\}\subset \H$ is a Bessel system for $\H$. Then, $\{f_j:j \in \N\}$ is complete in $\H$ if and only if there exists a Hilbert space $\H'$ containing $\H$ and a complete, $\omega$-linearly independent for $l^2(\N)$ system $\{e_j:j \in \N\}\subset \H'$ for $\H'$ such that
\[
\forall\; j\in \N, \quad P(e_j) = f_j.
\]

$\boldsymbol{e.}$ The collection $\{f_j:j \in \N\}\subset \H$ is a Bessel system for $\H$ with constant $B$ if and only if there exists a Hilbert space $\H'$ containing $\H$ and a collection $\{e_j:j \in \N\}\subset \H'$ satisfying the upper Riesz inequality with the same constant
and such that
\[
\forall\; j\in \N, \quad P(e_j) = f_j.
\]

\end{theorem}

\noindent{\bf Proof.}
$\boldsymbol{a.}$ This is the original result of Naimark.

$\boldsymbol{b.}$ This statement follows from the works of Casazza, Han, and Larson \cite{CHL}, Theorem 7.6, and Kashin and Kulikova \cite{KK}, Theorem 1.

$\boldsymbol{c.}$ This is the result of Terekhin, \cite{T}.

$\boldsymbol{d.}$
($\Longleftarrow$)
Let $\H \subset \H'$, let $\{e_j:j \in \N\}\subset \H'$ be a complete system for $\H'$, and let $f_j = P(e_j)$, $j\in \N$.
Then, the collection $\{f_j: j\in \N\} \subset \H$ is complete in $\H$.

($\Longrightarrow$)
As it was observed in \cite{KK}, it is enough to assume, without loss of generality, that $\H = l^2(\N)$.
Consider an infinite-dimensional matrix for which its $j$th column is defined to be the sequence of coefficients of $f_j \in \H$, $j\in\N$.
Let $v_i$, $i\in \N$ denote the rows of this matrix, i.e.,
\[
\forall\; i,j \in \N, \quad f_j(i) = v_i(j).
\]
With this definition, the assumption that $\{f_j: j \in \N\}$ forms a Bessel system for $\H$ implies that
\[
\forall\; i \in\N, \quad v_i \in l^2(\N),
\]
and, moreover, the system $\{v_i: i\in\N\}$ satisfies
\[
\forall\; c \in l^2(\N), \quad \Big\| \sum_{i\in\N} c(i) v_i \Big\|_2^2 \le B \|c\|_2^2.
\]

We first observe that $\{v_i:i\in\N\}$ is $\omega$-linearly independent for $l^2(\N)$-sequences in $\H$. Indeed, the fact that $\sum_{i\in\N} c(i) v_i =0$ for some $c \in l^2(\N)$ is equivalent to
\[
\forall\; j \in \N, \quad \langle c, f_j \rangle = 0.
\]
Thus, the completeness of $\{f_j:j\in\N\}$ implies that $c =0$. (This fact can also be derived from Proposition 4.2 in \cite{CKL}, once we observe that $\{v_i:i\in\N\}$ is the R-dual sequence, which corresponds to $\{f_j:j\in\N\}$.)

Let $V$ denote the closure in $l^2(\N)$ of the linear span of $\{v_i:i\in\N\}$. 
Let $V^\bot$ be the orthogonal complement of $V$ in $l^2(\N)$ and let $\{w_k: k \in K\}$ be an orthonormal basis of $V^\bot$, where $K$ denotes some set of indices (countable, finite, or empty, depending on the dimension of $V^\bot$), $K \cap \N = \emptyset$.
It is not difficult to see that the collection $\{v_i:i\in\N\}\cup \{w_k:k\in K\} \subset l^2(\N)$ is $\omega$-linearly independent for $l^2(\N)$-sequences in $\H$.

Consider now a matrix whose rows are the sequences of coefficients of vectors $\{v_i:i\in\N\}$ and $\{w_k:k\in K\}$, where the index set $L=\N\cup K$ of rows is endowed with some fixed order.

Let $e_j \in l^2(L) = \H'$ be the $j$th column of this matrix, $j\in\N$.
Clearly,
\[
\forall\; j\in \N, \quad P(e_j) = f_j.
\]

Since the row vectors $\{v_i:i\in\N\}\cup \{w_k:k\in K\} \subset l^2(\N)$ are $\omega$-linearly independent for $l^2(\N)$-sequences in $\H$, it follows that the collection $\{e_j:j \in \N\}$ is complete in $\H' = l^2(L)$.

In order to show that $\{e_j: j\in\N\}$ is $\omega$-linearly independent for $l^2(\N)$, let $c \in l^2(\N)$ be such that
\[
\sum_{j\in\N} c(j) e_j =0.
\]
This implies that
\begin{equation}
\label{e1}
\forall\; i \in \N, \quad \langle c, v_i \rangle =0
\end{equation}
and
\begin{equation}
\label{e2}
\forall\; k \in K, \quad \langle c, w_k \rangle =0.
\end{equation}
Write $c = c_V + c_{V^\bot}$, where $c_V \in V$ and $c_{V^\bot} \in V^\bot$. Then, since $\{w_k:k\in K\}$ is an orthonormal basis in $V^\bot$, (\ref{e2}) implies that $c_{V^\bot} =0$. Then, (\ref{e1}) may be rewritten as:
\[
\forall\; i \in \N, \quad \langle c_V, v_i \rangle =0.
\]
By the definition of the space $V$, $\{v_i:i\in \N\}$ is complete in $V$, and so $c_V =0$, as well.

$\boldsymbol{e.}$
($\Longleftarrow$)
Let $\H \subset \H'$, let $\{e_j:j \in \N\}\subset \H'$ satisfy the second inequality in (\ref{upRiesz}) for $\H'$, and let $f_j = P(e_j)$, $j\in \N$.
In particular, it follows that $\{e_j:j \in \N\}$ is a Bessel system for $\H'$. Thus, for each $h \in \H$ we have
\[
\sum_{j\in \N} |\langle h , P(e_j) \rangle|^2 = \sum_{j\in \N} |\langle P^*(h) , e_j \rangle|^2 \le B\|P^*(h)\|_2^2 \le B \|h\|_2^2.
\]

($\Longrightarrow$)
Again, we assume without loss of generality that $\H = l^2(\N)$.
As we have shown in the proof of part {\it d}, the assumption that $\{f_j:j\in \N\}$ is a Bessel system implies that
\begin{equation}
\label{e3}
\forall\; c \in l^2(\N), \quad \Big\| \sum_{i\in \N} c(i) v_i \Big\|_2^2 \le B \|c\|_2^2.
\end{equation}
As before, we consider the decomposition $l^2(\N) = V \oplus V^\bot$, we choose an orthogonal basis $\{w_k:k\in K\}$ for $V^\bot$, such that $\|w_k\| \le B$, $k \in K$, and we define the collection $\{e_l: l \in L\} \subset \H' = l^2(L)$. Inequality (\ref{e3}) implies that $\{v_i:i \in \N\}$ is a Bessel system with constant $B$; so is $\{w_k:k\in K\}$.
Therefore, we obtain the desired estimate:
\[
\begin{aligned}
\Big\| \sum_{l\in L} c(l) e_l \Big\|_2^2 &= \sum_{l\in L} \sum_{i\in \N} |c(l) e_l(i)|^2\\
&= \sum_{i\in \N} |\langle c_V, v_i\rangle|^2 + \sum_{k\in K} |\langle c_{V^\bot}, w_k\rangle|^2\\
&\le B\|c\|_2^2.
\end{aligned}
\]


\hfill \qed


\end{document}